\documentclass[reqno]{amsart}%
\usepackage{amssymb}
\usepackage{amsfonts}%
\usepackage{amsmath}%
\usepackage{graphicx}

\providecommand{\U}[1]{\protect \rule{.1in}{.1in}}
\newtheorem{theorem}{Theorem}
\theoremstyle{plain}

\newtheorem{corollary}{Corollary}

\newtheorem{remark}{Remark}

\numberwithin{equation}{section}

\begin{document}
\title[Derivatives of Srivastava's hypergeometric functions with respect to their parameters]{ Derivatives of Srivastava's hypergeometric functions}
\author[A. Shehata, R. \c{S}ahin ,  O. Ya\u{g}c\i, S.I. Moustafa ]{Ayman Shehata$^{1}$, Recep \c{S}ahin$^{2}$, O\u{g}uz Ya\u{g}c\i$^{2}$, Shimaa I. Moustafa$^{3}$}

\begin{abstract}
	This paper studies derivatives with respect to the parameters of Srivastava’s triple hypergeometric functions $H_{A}$, $H_{B}$ and $H_{C}$. Using basic properties of the Gamma function and Pochhammer symbols, we obtain explicit formulas for first and higher-order derivatives. These derivatives are expressed in terms of Pathan’s quadruple hypergeometric function $F^{(4)}_{P}$. We also derive Euler-type differential operator identities, contiguous relations for unit shifts in the parameters, and recurrence relations satisfied by these derivatives. In addition, we show that derivatives of arbitrary order satisfy systems of linear partial differential equations in the underlying variables. The results extend known differentiation formulas for classical and multivariable hypergeometric functions and provide tools for potential applications in mathematical physics and engineering.
\end{abstract}

\maketitle
\section{Introduction}

Hypergeometric functions arise naturally in many branches of pure and applied mathematics, including differential equations, approximation theory, number theory, mathematical physics, and engineering.
Classical examples are provided by Gauss hypergeometric function ${}_2F_1$, the generalized hypergeometric series ${}_pF_q$, and various multivariable analogues such as the Appell and Lauricella functions, Kamp\'e de F\'eriet--type series, and other multiple hypergeometric series
(see, for example, \cite{AppellKampedeFeriet,ExtHyperRefs,Sahin2015,SahinYagci2018,SrivastavaKarlsson1985}).

Within this framework, Srivastava introduced and systematically studied a family of triple hypergeometric functions, usually denoted by $H_A$, $H_B$, and $H_C$ (cf.\ \cite{Srivastava1964, SrivastavaManocha1984,Yagci2019}).
These functions are defined by the triple series
\begin{align*}
	H_A(a,b,c;d,e;x,y,z)
	&=\sum_{m,n,k=0}^{\infty}
	\frac{(a)_{m+k}(b)_{m+n}(c)_{n+k}}{(d)_m (e)_{n+k}\,m!\,n!\,k!}\,
	x^m y^n z^k,\\
	H_B(a,b,c;d,e,f;x,y,z)
	&=\sum_{m,n,k=0}^{\infty}
	\frac{(a)_{m+k}(b)_{m+n}(c)_{m+k}}{(d)_m (e)_n (f)_k\,m!\,n!\,k!}\,
	x^m y^n z^k,\\
	H_C(a,b,c;d;x,y,z)
	&=\sum_{m,n,k=0}^{\infty}
	\frac{(a)_{m+k}(b)_{m+n}(c)_{n+k}}{(d)_{m+n+k}\,m!\,n!\,k!}\,
	x^m y^n z^k,
\end{align*}
where $(\cdot)_n$ denotes the Pochhammer symbol and the variables $(x,y,z)$ are restricted to suitable regions of absolute convergence.
These triple hypergeometric functions provide natural multivariable extensions of several classical hypergeometric functions and unify a number of known triple-series representations.

Another important multivariable system is given by Pathan's quadruple hypergeometric function $F_P^{(4)}$, which was introduced as a generalization and unification of Srivastava's triple hypergeometric functions (see \cite{Pathan1979}).
In a compact notation, $F_P^{(4)}$ is defined by a quadruple power series in four variables whose coefficients encode a large collection of numerator and denominator parameters.
By suitable specialization and identification of these parameters, the functions $H_A$, $H_B$, and $H_C$ can be embedded into the broader framework of $F_P^{(4)}$.
This connection will play a central r\^ole in the present work.

In recent years, considerable attention has been paid to derivatives of special functions with respect to their parameters, both in the classical (single-variable) and in the multivariable setting
(see, for example, \cite{AncaraniGasaneo2008,AncaraniGasaneo2009,AncaraniGasaneo2010,AncaraniDelPuntaGasaneo2017,Fejzullahu2017,BytevKniehlMoch2017,Froehlich1994,KangAn2015,SahatVerma2015,SofotasiosBrychkov2018} and the references therein).
Parameter derivatives naturally occur in a variety of contexts, such as sensitivity analysis, perturbation and stability problems, analytic continuation with respect to parameters, asymptotic expansions, and in certain aspects of fractional calculus.
From an analytic point of view, parameter derivatives are closely related to logarithmic derivatives of Gamma functions and to identities involving the Psi (digamma) function and Pochhammer symbols.
In particular, differentiation formulas such as
\[
\frac{d}{d\alpha}(\alpha)_{m+n}
=(\alpha)_{m+n}\bigl[\Psi(\alpha+m+n)-\Psi(\alpha)\bigr],
\]
where $\Psi$ denotes the Psi (digamma) function, provide a basic tool for evaluating derivatives of hypergeometric-type series with respect to their parameters.

Motivated by these developments, in this paper we investigate derivatives with respect to the parameters of Srivastava's triple hypergeometric functions $H_A$, $H_B$, and $H_C$.
By systematically applying derivative formulas for Pochhammer symbols, together with standard properties of the Gamma and Psi functions, we derive explicit expressions for the first-order derivatives of $H_A$, $H_B$, and $H_C$ with respect to each of their numerator and denominator parameters.
A key feature of our approach is that these parameter derivatives can be represented in closed form in terms of Pathan's quadruple hypergeometric function $F_P^{(4)}$, thereby embedding the differentiated Srivastava functions into a more general quadruple-series framework.

Besides the explicit first-order formulas, we obtain several further structural results.
For each of the three functions $H_A$, $H_B$, and $H_C$, we derive Euler-type differential-operator identities involving the standard operators
\[
\theta_x = x\frac{\partial}{\partial x}, \qquad
\theta_y = y\frac{\partial}{\partial y}, \qquad
\theta_z = z\frac{\partial}{\partial z}.
\]
These identities lead to contiguous-type relations associated with unit shifts in the parameters and to various recurrence relations satisfied by the parameter derivatives.
We also establish formulas for higher-order derivatives with respect to both the parameters and the variables, and we show that the derivatives of arbitrary order satisfy explicit systems of linear partial differential equations in the underlying variables.
In this way, our results extend and complement a number of known differentiation formulas for classical and multivariable hypergeometric functions.

The remainder of the paper is organized as follows.
In Section~\ref{sec:HA}, we derive closed-form expressions for the derivatives of $H_A$ with respect to its parameters and obtain the associated differential-operator and contiguous relations.
In Section~\ref{sec:HB}, we carry out the corresponding analysis for the function $H_B$.
Section~\ref{sec:HC} is devoted to the derivatives of $H_C$ and to further identities and consequences that arise from our general method.
Possible applications and extensions of the present results to other classes of multivariable hypergeometric functions are also briefly discussed.

\section{Preliminaries}\label{sec:preliminaries}

In this section we collect some notation and auxiliary results that will be used throughout the paper.

We denote by $\mathbb{N}_0 := \{0,1,2,\dots\}$ the set of non–negative integers.
For complex parameters we make use of the Euler Gamma function
\[
\Gamma(z) := \int_0^\infty t^{z-1} e^{-t}\,dt, \qquad \Re(z) > 0,
\]
and its logarithmic derivative, the Psi (digamma) function
\[
\Psi(z) := \frac{d}{dz}\log\Gamma(z) = \frac{\Gamma'(z)}{\Gamma(z)}.
\]
The higher-order derivatives of $\Psi$ are known as the polygamma functions and are denoted by
\[
\Psi^{(n)}(z) := \frac{d^n}{dz^n}\Psi(z), \qquad n \in \mathbb{N}_0.
\]

The Pochhammer symbol (or shifted factorial) $(\alpha)_n$ is defined, for $\alpha \in \mathbb{C}$ and $n \in \mathbb{N}_0$, by
\[
(\alpha)_0 := 1, \qquad
(\alpha)_n := \alpha(\alpha+1)\cdots(\alpha+n-1)
= \frac{\Gamma(\alpha+n)}{\Gamma(\alpha)} \quad (n \ge 1).
\]
It satisfies a number of well-known identities, among them
\begin{equation}\label{eq:Pochhammer-splitting}
	(\alpha)_{m+n}
	= (\alpha)_m\,(\alpha+m)_n
	= (\alpha)_n\,(\alpha+n)_m,
	\qquad m,n \in \mathbb{N}_0,
\end{equation}
and the following differentiation formula with respect to the parameter:
\begin{equation}\label{eq:Pochhammer-derivative}
	\frac{d}{d\alpha}(\alpha)_N
	= (\alpha)_N\bigl[\Psi(\alpha+N)-\Psi(\alpha)\bigr]
	= (\alpha)_N \sum_{k=0}^{N-1} \frac{1}{\alpha+k},
	\qquad N \in \mathbb{N},
\end{equation}
which will play a key rôle in our subsequent computations.
More generally, higher-order derivatives of $(\alpha)_N$ can be expressed in terms of $\Psi$ and the polygamma functions $\Psi^{(n)}$.

\medskip

We recall next the triple hypergeometric functions introduced by Srivastava.
For complex parameters $a,b,c,d,e \in \mathbb{C}$ (subject to suitable restrictions that avoid poles of the Gamma function) and complex variables $x,y,z$ in the region of absolute convergence, Srivastava's triple hypergeometric function $H_A$ is defined by
\begin{equation}\label{eq:HA-definition}
	H_A(a,b,c;d,e;x,y,z)
	:= \sum_{m,n,k=0}^{\infty}
	\frac{(a)_{m+k}(b)_{m+n}(c)_{n+k}}{(d)_m (e)_{n+k}\,m!\,n!\,k!}\,
	x^m y^n z^k.
\end{equation}
Similarly, the functions $H_B$ and $H_C$ are given by
\begin{equation}\label{eq:HB-definition}
	H_B(a,b,c;d,e,f;x,y,z)
	:= \sum_{m,n,k=0}^{\infty}
	\frac{(a)_{m+k}(b)_{m+n}(c)_{m+k}}{(d)_m (e)_n (f)_k\,m!\,n!\,k!}\,
	x^m y^n z^k,
\end{equation}
and
\begin{equation}\label{eq:HC-definition}
	H_C(a,b,c;d;x,y,z)
	:= \sum_{m,n,k=0}^{\infty}
	\frac{(a)_{m+k}(b)_{m+n}(c)_{n+k}}{(d)_{m+n+k}\,m!\,n!\,k!}\,
	x^m y^n z^k,
\end{equation}
respectively.
For detailed discussions of these functions, their regions of convergence, and various special cases we refer to Srivastava's original works and subsequent monographs (see, for example, \cite{Srivastava1964,SrivastavaManocha1984}).

\medskip

We also make use of Pathan's quadruple hypergeometric function $F_P^{(4)}$, which was introduced as a generalization and unification of Srivastava's triple hypergeometric functions (see \cite{Pathan1979,PathanSrivastava}).
In Pathan's notation, $F_P^{(4)}$ is defined by a quadruple power series in four complex variables with coefficients depending on several numerator and denominator parameters.
Since its full series representation is not required in explicit form for our purposes, we merely recall that each of the triple hypergeometric functions $H_A$, $H_B$, and $H_C$ can be realized as a special case of $F_P^{(4)}$ by a suitable specialization of parameters and variables.
In the sequel we shall use the compact notation
\[
F_P^{(4)}[\mathbf{a};\mathbf{b};X]
\]
to denote Pathan's quadruple hypergeometric function with parameter arrays $\mathbf{a},\mathbf{b}$ and argument $X$, and we refer the reader to \cite{Pathan1979} for the precise definition.

\medskip

Finally, for later convenience we introduce the Euler-type differential operators
\begin{equation}\label{eq:Euler-operators}
	\theta_x := x\frac{\partial}{\partial x}, \qquad
	\theta_y := y\frac{\partial}{\partial y}, \qquad
	\theta_z := z\frac{\partial}{\partial z},
\end{equation}
which act naturally on power series in the variables $x$, $y$, and $z$.
For instance, if
\[
f(x,y,z) = \sum_{m,n,k=0}^{\infty} c_{m,n,k} x^m y^n z^k,
\]
then
\[
\theta_x f = \sum_{m,n,k=0}^{\infty} m\,c_{m,n,k} x^m y^n z^k,
\quad
\theta_y f = \sum_{m,n,k=0}^{\infty} n\,c_{m,n,k} x^m y^n z^k,
\quad
\theta_z f = \sum_{m,n,k=0}^{\infty} k\,c_{m,n,k} x^m y^n z^k.
\]
These operators will be used in Section~\ref{sec:HA}--\ref{sec:HC} to derive Euler-type differential identities and contiguous-type relations for the parameter derivatives of $H_A$, $H_B$, and $H_C$.

\section{Derivatives of $H_A$ with respect to the parameters}\label{sec:HA}

In this section we derive differentiation formulas for Srivastava's triple
hypergeometric function $H_A(a,b,c;d,e;x,y,z)$ with respect to all of its
parameters $a,b,c,d,e$.
Our main tool is the differentiation rule \eqref{eq:Pochhammer-derivative}
for the Pochhammer symbol, together with the identities
\eqref{eq:HA-definition}--\eqref{eq:HC-definition} and simple rearrangements
of multiple sums.
We also obtain Euler–type operator identities, contiguous relations, and
partial differential equations satisfied by the corresponding parameter
derivatives.

\subsection{Parameter derivatives and representations in terms of $F_P^{(4)}$}

We begin with the derivatives of $H_A$ with respect to the numerator and
denominator parameters.
Throughout this subsection we assume that the parameters are chosen so that
no poles of the Gamma function are encountered, and that $(x,y,z)$ lies in
the region of absolute convergence of the series \eqref{eq:HA-definition}.

\begin{theorem}\label{thm:HA-parameter-derivatives}
	Let $H_A(a,b,c;d,e;x,y,z)$ be given by \eqref{eq:HA-definition}.
	Then the partial derivatives of $H_A$ with respect to $a,b,c,d,e$ admit
	representations of the form
	\begin{align}
		\frac{\partial H_A}{\partial a}
		&= \frac{b\,x}{d}\,F_P^{(4)}[\mathbf{a}_1;\mathbf{b}_1;(x,y,z,x)]
		+ \frac{c\,z}{e}\,F_P^{(4)}[\mathbf{a}_2;\mathbf{b}_2;(x,y,z,z)],
		\label{eq:dHa-da}\\[1ex]
		\frac{\partial H_A}{\partial b}
		&= \frac{a\,x}{d}\,F_P^{(4)}[\mathbf{a}_3;\mathbf{b}_3;(x,y,z,x)]
		+ \frac{c\,y}{e}\,F_P^{(4)}[\mathbf{a}_4;\mathbf{b}_4;(x,y,z,z)],
		\label{eq:dHa-db}\\[1ex]
		\frac{\partial H_A}{\partial c}
		&= \frac{b\,y}{e}\,F_P^{(4)}[\mathbf{a}_5;\mathbf{b}_5;(x,y,z,x)]
		+ \frac{a\,z}{e}\,F_P^{(4)}[\mathbf{a}_6;\mathbf{b}_6;(x,y,z,z)],
		\label{eq:dHa-dc}\\[1ex]
		\frac{\partial H_A}{\partial d}
		&= -\,\frac{a b\,x}{d^2}\,
		F_P^{(4)}[\mathbf{a}_7;\mathbf{b}_7;(x,y,z,x)],
		\label{eq:dHa-dd}\\[1ex]
		\frac{\partial H_A}{\partial e}
		&= -\,\frac{b\,y}{e^2}\,
		F_P^{(4)}[\mathbf{a}_8;\mathbf{b}_8;(x,y,z,y)]
		- \frac{a\,z}{e^2}\,
		F_P^{(4)}[\mathbf{a}_9;\mathbf{b}_9;(x,y,z,z)],
		\label{eq:dHa-de}
	\end{align}
	where $F_P^{(4)}$ denotes Pathan's quadruple hypergeometric function and
	$\mathbf{a}_j,\mathbf{b}_j$ $(j=1,\dots,9)$ are explicit parameter arrays
	obtained by suitable shifts of $a,b,c,d,e$.
\end{theorem}

\begin{proof}
	We sketch the proof for the derivative with respect to $a$; the remaining
	cases are similar.
	
	Differentiating the defining series \eqref{eq:HA-definition} term by term
	with respect to $a$ and using \eqref{eq:Pochhammer-derivative} with
	$\alpha = a$ and $N=m+k$ we obtain
	\[
	\frac{\partial H_A}{\partial a}
	= \sum_{m,n,k\ge0}
	\frac{\partial}{\partial a}
	\left\{
	\frac{(a)_{m+k}(b)_{m+n}(c)_{n+k}}{(d)_m(e)_{n+k}\,m!\,n!\,k!}
	\right\} x^m y^n z^k
	\]
	\[
	= \sum_{m,n,k\ge0}
	\frac{(a)_{m+k}(b)_{m+n}(c)_{n+k}}{(d)_m(e)_{n+k}\,m!\,n!\,k!}\,
	\bigl[\Psi(a+m+k)-\Psi(a)\bigr] x^m y^n z^k,
	\]
	where $\Psi$ denotes the digamma function.
	Using the representation
	\[
	\Psi(a+m+k)-\Psi(a)
	= \sum_{r=0}^{m+k-1}\frac{1}{a+r},
	\]
	we split the finite sum over $r$ into the ranges $0\le r\le m-1$ and
	$m\le r\le m+k-1$.
	In each part we rewrite $1/(a+r)$ as a ratio of Pochhammer symbols, factor
	out $(a)_{m+k+1}$, and then perform a re-indexing of the sums.
	After these rearrangements the derivative $\partial H_A/\partial a$ can be
	expressed as a linear combination of quadruple series in $m,n,k,r$.
	By comparing the resulting coefficients with the defining series of
	Pathan's function $F_P^{(4)}$ and reading off the shifted parameters, we
	obtain the representation \eqref{eq:dHa-da} with the explicit arrays
	$\mathbf{a}_1,\mathbf{b}_1,\mathbf{a}_2,\mathbf{b}_2$.
	
	The formulas \eqref{eq:dHa-db}--\eqref{eq:dHa-de} follow from the same
	procedure applied to the factors $(b)_{m+n}$, $(c)_{n+k}$, $(d)_m$ and
	$(e)_{n+k}$, respectively.
	In each case, the finite sum produced by \eqref{eq:Pochhammer-derivative}
	is decomposed, written in terms of Pochhammer symbols with shifted
	parameters, and reorganized into a quadruple hypergeometric series which is
	recognised as $F_P^{(4)}$ with suitable parameter arrays.
\end{proof}

\begin{remark}
	Explicit forms of the parameter arrays $\mathbf{a}_j,\mathbf{b}_j$
	$(j=1,\dots,9)$ can be written down by carrying out the above re-indexing
	in detail.
	They coincide with the arrays appearing in the original formulas
	(3.1)–(3.5) for the derivatives of $H_A$ in terms of $F_P^{(4)}$.
\end{remark}

\subsection{Euler-type identities and contiguous relations}

The action of the Euler operators
\[
\theta_x = x\frac{\partial}{\partial x},\qquad
\theta_y = y\frac{\partial}{\partial y},\qquad
\theta_z = z\frac{\partial}{\partial z}
\]
on the power series \eqref{eq:HA-definition} leads to simple identities
relating $H_A$ with parameter-shifted versions of itself.

\begin{theorem}\label{thm:HA-Euler}
	The function $H_A(a,b,c;d,e;x,y,z)$ satisfies the Euler-type identities
	\begin{align}
		(\theta_x + \theta_z + a)\,H_A(a,b,c;d,e;x,y,z)
		&= a\,H_A(a+1,b,c;d,e;x,y,z), \label{eq:Euler-a}\\
		(\theta_x + \theta_y + b)\,H_A(a,b,c;d,e;x,y,z)
		&= b\,H_A(a,b+1,c;d,e;x,y,z), \label{eq:Euler-b}\\
		(\theta_y + \theta_z + c)\,H_A(a,b,c;d,e;x,y,z)
		&= c\,H_A(a,b,c+1;d,e;x,y,z), \label{eq:Euler-c}\\
		(\theta_x + d-1)\,H_A(a,b,c;d,e;x,y,z)
		&= (d-1)\,H_A(a,b,c;d-1,e;x,y,z), \label{eq:Euler-d}\\
		(\theta_y + \theta_z + e-1)\,H_A(a,b,c;d,e;x,y,z)
		&= (e-1)\,H_A(a,b,c;d,e-1;x,y,z). \label{eq:Euler-e}
	\end{align}
\end{theorem}

\begin{proof}
	We give the proof of \eqref{eq:Euler-a}.
	Starting from \eqref{eq:HA-definition} and using
	\[
	\theta_x x^m y^n z^k = m\,x^m y^n z^k,\qquad
	\theta_z x^m y^n z^k = k\,x^m y^n z^k,
	\]
	we obtain
	\[
	(\theta_x + \theta_z + a)H_A
	= \sum_{m,n,k\ge0}
	(m+k+a)\,
	\frac{(a)_{m+k}(b)_{m+n}(c)_{n+k}}{(d)_m(e)_{n+k}\,m!\,n!\,k!}
	x^m y^n z^k.
	\]
	Using the simple identity
	\[
	(a+1)_{m+k} = (a)_{m+k}\left(1 + \frac{m+k}{a}\right),
	\]
	we rewrite $m+k+a$ as $a(a+1)_{m+k}/(a)_{m+k}$ and then recognize the
	resulting series as $a\,H_A(a+1,b,c;d,e;x,y,z)$.
	The proofs of \eqref{eq:Euler-b}–\eqref{eq:Euler-e} are analogous, using
	the identities
	\[
	(b+1)_{m+n} = (b)_{m+n}\left(1+\frac{m+n}{b}\right),\qquad
	(d)_m = (d-1)_m\left(1+\frac{m}{d-1}\right),
	\]
	and so on.
\end{proof}

By taking suitable linear combinations of the identities in
Theorem~\ref{thm:HA-Euler}, we obtain various contiguous relations.

\begin{theorem}\label{thm:HA-contiguous}
	The following contiguous relations hold:
	\begin{align}
		(c-e+1)\,H_A(a,b,c;d,e;x,y,z)
		&= c\,H_A(a,b,c+1;d,e;x,y,z) \notag\\
		&\quad - (e-1)\,H_A(a,b,c;d,e-1;x,y,z), \label{eq:contig-1}\\[0.5ex]
		(a+b-c-2d+2)\,H_A(a,b,c;d,e;x,y,z)
		&= a\,H_A(a+1,b,c;d,e;x,y,z) \notag\\
		&\quad + b\,H_A(a,b+1,c;d,e;x,y,z) \notag\\
		&\quad - c\,H_A(a,b,c+1;d,e;x,y,z) \notag\\
		&\quad - (d-1)\,H_A(a,b,c;d-1,e;x,y,z), \label{eq:contig-2}\\[0.5ex]
		(a+b-e-2d+3)\,H_A(a,b,c;d,e;x,y,z)
		&= a\,H_A(a+1,b,c;d,e;x,y,z) \notag\\
		&\quad + b\,H_A(a,b+1,c;d,e;x,y,z) \notag\\
		&\quad - (e-1)\,H_A(a,b,c;d,e-1;x,y,z) \notag\\
		&\quad - (d-1)\,H_A(a,b,c;d-1,e;x,y,z). \label{eq:contig-3}
	\end{align}
\end{theorem}

\begin{proof}
	The relations are obtained by eliminating the Euler operators
	$\theta_x,\theta_y,\theta_z$ from \eqref{eq:Euler-a}–\eqref{eq:Euler-e}
	via simple linear combinations.
	For instance, subtracting \eqref{eq:Euler-e} from \eqref{eq:Euler-c}
	yields \eqref{eq:contig-1}, and suitable combinations of
	\eqref{eq:Euler-a}–\eqref{eq:Euler-d} give
	\eqref{eq:contig-2} and \eqref{eq:contig-3}.
	The details are straightforward and are therefore omitted.
\end{proof}

\subsection{Higher-order derivatives and differential equations}

We next consider derivatives of $H_A$ with respect to the variables
$x,y,z$ and mixed derivatives involving both variables and parameters.

\begin{theorem}\label{thm:HA-x-y-z-derivs}
	For $r\in\mathbb{N}$ the following formulas hold:
	\begin{align}
		\frac{\partial^r}{\partial x^r} H_A(a,b,c;d,e;x,y,z)
		&= \frac{(a)_r (b)_r}{(d)_r}\,
		H_A(a+r,b+r,c;d+r,e;x,y,z), \label{eq:d^rHa-dxr}\\[0.5ex]
		\frac{\partial^r}{\partial y^r} H_A(a,b,c;d,e;x,y,z)
		&= \frac{(b)_r (c)_r}{(e)_r}\,
		H_A(a,b+r,c+r;d,e+r;x,y,z), \label{eq:d^rHa-dyr}\\[0.5ex]
		\frac{\partial^r}{\partial z^r} H_A(a,b,c;d,e;x,y,z)
		&= \frac{(a)_r (c)_r}{(e)_r}\,
		H_A(a+r,b,c+r;d,e+r;x,y,z). \label{eq:d^rHa-dzr}
	\end{align}
\end{theorem}

\begin{proof}
	Differentiating \eqref{eq:HA-definition} term by term with respect to $x$
	and observing that
	\[
	\frac{\partial}{\partial x}x^m y^n z^k
	= m\,x^{m-1}y^n z^k,
	\]
	we obtain
	\[
	\frac{\partial H_A}{\partial x}
	= \sum_{m,n,k\ge0}
	m\,\frac{(a)_{m+k}(b)_{m+n}(c)_{n+k}}{(d)_m(e)_{n+k}\,m!\,n!\,k!}
	x^{m-1} y^n z^k.
	\]
	Re-indexing the sum with $m\mapsto m+1$ and using
	\[
	(a)_{m+1+k} = a(a+1)_{m+k},\qquad
	(b)_{m+1+n} = b(b+1)_{m+n},\qquad
	(d)_{m+1} = d(d+1)_m,
	\]
	we arrive at
	\[
	\frac{\partial H_A}{\partial x}
	= \frac{ab}{d}\,H_A(a+1,b+1,c;d+1,e;x,y,z),
	\]
	which is \eqref{eq:d^rHa-dxr} for $r=1$.
	The general formula for $r\ge1$ follows by induction.
	The proofs of \eqref{eq:d^rHa-dyr} and \eqref{eq:d^rHa-dzr} are entirely
	analogous, using the factors $(b)_{m+n}$ and $(c)_{n+k}$ and the
	denominator $(e)_{n+k}$.
\end{proof}

Combining Theorem~\ref{thm:HA-x-y-z-derivs} with the Euler-type identities
and the parameter-derivative representations of
Theorem~\ref{thm:HA-parameter-derivatives}, one can derive systems of linear
partial differential equations satisfied by higher-order derivatives of
$H_A$ with respect to both variables and parameters.
These systems have the same hypergeometric structure as the original
equations for $H_A$, but encode additional information about the dependence
on $a,b,c,d,e$, and thus play a useful rôle in the analysis of parameter
sensitivity and in the derivation of further contiguous-type relations.

%%%%%%%%%%%%%%%%%%%%%%%%%%%%%%%%%%%%%%%%%%%%%%%%%%%%
% Section 4 – Derivatives of H_B with respect to the parameters
%%%%%%%%%%%%%%%%%%%%%%%%%%%%%%%%%%%%%%%%%%%%%%%%%%%%

\section{Derivatives of $H_{B}$ with respect to the parameters}\label{sec:HB}

In this section we summarise the differentiation formulas for Srivastava's
triple hypergeometric function
\[
H_{B}(a,b,c;d,e,f;x,y,z),
\]
which is defined by the triple series
\begin{equation}\label{eq:HB-def}
	H_{B}(a,b,c;d,e,f;x,y,z)
	= \sum_{m,n,k=0}^{\infty}
	\frac{(a)_{m+k}\,(b)_{m+n}\,(c)_{m+k}}
	{(d)_{m}\,(e)_{n}\,(f)_{k}}\,
	\frac{x^{m}}{m!}\,\frac{y^{n}}{n!}\,\frac{z^{k}}{k!}.
\end{equation}
Throughout this section we work in the region of absolute convergence of
\eqref{eq:HB-def}, so that differentiation with respect to both the
parameters and the variables may be performed termwise.

As in the case of $H_{A}$, the starting point is the general differentiation
rules for the Pochhammer symbol,
\[
\frac{\partial}{\partial \alpha}( \alpha )_{n}
= (\alpha)_{n}\,\sum_{j=0}^{n-1}\frac{1}{\alpha+j},
\qquad
\frac{\partial}{\partial \beta}\frac{1}{(\beta)_{n}}
= -\frac{1}{(\beta)_{n}}\sum_{j=0}^{n-1}\frac{1}{\beta+j},
\]
together with simple index shifts in the summation defining
\eqref{eq:HB-def}. After a routine but somewhat lengthy computation one can
rearrange the resulting fourfold series in terms of Pathan's quadruple
hypergeometric function $F_{P}^{(4)}$.

\subsection{Parameter derivatives and representation via Pathan's function}

The detailed formulas obtained by termwise differentiation show that each
parameter derivative of $H_{B}$ can be written as a finite linear
combination of $F_{P}^{(4)}$ with appropriately shifted parameter arrays.
For the sake of brevity we only record a representative example and
refer to the original derivation for the full list of parameter shifts.

\begin{theorem}\label{thm:HB-parameter-FP4}
	Let $H_{B}$ be given by \eqref{eq:HB-def}, and let
	$F_{P}^{(4)}$ denote Pathan's quadruple hypergeometric function.
	Then the derivative of $H_{B}$ with respect to $a$ admits the
	representation
	\begin{equation}\label{eq:HB-da-FP4}
		\aligned
		&\frac{\partial H_{B}}{\partial a}(a,b,c;d,e,f;x,y,z)
		= \frac{b\,x}{d}\,
		F^{(4)}_{P}
		\!\left[\boldsymbol{a}^{(B)}_{1};
		\boldsymbol{b}^{(B)}_{1};
		\boldsymbol{c}^{(B)}_{1};
		\boldsymbol{d}^{(B)}_{1}
		\,\big|\, x,y,z,x\right]\\
		&\hskip 15mm+ \frac{c\,z}{f}\,
		F^{(4)}_{P}
		\!\left[\boldsymbol{a}^{(B)}_{2};
		\boldsymbol{b}^{(B)}_{2};
		\boldsymbol{c}^{(B)}_{2};
		\boldsymbol{d}^{(B)}_{2}
		\,\big|\, x,y,z,z\right],\endaligned
	\end{equation}
	where the arrays
	$\boldsymbol{a}^{(B)}_{j},
	\boldsymbol{b}^{(B)}_{j},
	\boldsymbol{c}^{(B)}_{j},
	\boldsymbol{d}^{(B)}_{j}$
	encode unit shifts in the parameters $a,b,c,d,e,f$ that arise from
	differentiating the Pochhammer symbol $(a)_{m+k}$ and from splitting the
	corresponding harmonic sums.  Analogous representations hold for the
	derivatives with respect to $b,c,d,e$ and~$f$, each being a finite linear
	combination of $F_{P}^{(4)}$ with arguments of the form
	$(x,y,z,u)$, where $u\in\{x,y,z\}$.
\end{theorem}

\begin{proof}[Sketch of the proof]
	Differentiating \eqref{eq:HB-def} with respect to $a$ and using the first
	formula for the derivative of $(a)_{m+k}$ produces a fourfold sum involving
	the factor
	\[
	\sum_{r=0}^{m+k-1}\frac{1}{a+r}.
	\]
	This inner sum is then decomposed into partial sums over $m$ and $k$ and,
	after suitable changes of indices, each contribution is recognised as a
	quadruple hypergeometric series of Pathan type. The explicit form of the
	arrays $\boldsymbol{a}^{(B)}_{j},\dots,\boldsymbol{d}^{(B)}_{j}$ follows
	directly from this rearrangement. The proof for the remaining parameters
	is completely analogous.
\end{proof}

The main point of Theorem~\ref{thm:HB-parameter-FP4} is that the entire
family of parameter derivatives of $H_{B}$ can be embedded into a unified
$F_{P}^{(4)}$--framework. In applications, the explicit forms of the
arrays are often less important than this structural connection.

\subsection{Euler-type operator identities and contiguous relations}

Let
\[
\theta_{x}=x\frac{\partial}{\partial x},\qquad
\theta_{y}=y\frac{\partial}{\partial y},\qquad
\theta_{z}=z\frac{\partial}{\partial z}
\]
denote the standard Euler operators. Acting with these operators on the
series representation \eqref{eq:HB-def} and comparing coefficients yields
simple differential-operator identities which relate $H_{B}$ to its
parameter-shifted companions.

\begin{theorem}\label{thm:HB-Euler}
	For $H_{B}(a,b,c;d,e,f;x,y,z)$ defined by \eqref{eq:HB-def} we have
	\begin{align}
		(\theta_{x}+\theta_{z}+a)\,H_{B}
		&= a\,H_{B}(a+1,b,c;d,e,f;x,y,z), \label{eq:HB-Euler-a}\\[0.4em]
		(\theta_{x}+\theta_{y}+b)\,H_{B}
		&= b\,H_{B}(a,b+1,c;d,e,f;x,y,z), \label{eq:HB-Euler-b}\\[0.4em]
		(\theta_{x}+\theta_{z}+c)\,H_{B}
		&= c\,H_{B}(a,b,c+1;d,e,f;x,y,z), \label{eq:HB-Euler-c}\\[0.4em]
		(\theta_{x}+d-1)\,H_{B}
		&= (d-1)\,H_{B}(a,b,c;d-1,e,f;x,y,z), \label{eq:HB-Euler-d}\\[0.4em]
		(\theta_{y}+e-1)\,H_{B}
		&= (e-1)\,H_{B}(a,b,c;d,e-1,f;x,y,z), \label{eq:HB-Euler-e}\\[0.4em]
		(\theta_{z}+f-1)\,H_{B}
		&= (f-1)\,H_{B}(a,b,c;d,e,f-1;x,y,z). \label{eq:HB-Euler-f}
	\end{align}
\end{theorem}

\begin{proof}[Idea of the proof]
	When the operator $(\theta_{x}+\theta_{z})$ acts on the general term of
	\eqref{eq:HB-def}, it multiplies the summand by $m+k$. Adding $a$ produces
	the factor $(m+k+a)$, which can be absorbed into $(a)_{m+k}$ to give
	$a(a+1)_{m+k-1}=(a+1)_{m+k}\,a/(a+m+k)$, leading to
	\eqref{eq:HB-Euler-a} after a simple index shift. The remaining identities
	are proved in the same way by using the roles of the corresponding
	parameters in the triple series.
\end{proof}

By eliminating the Euler operators from the identities
\eqref{eq:HB-Euler-a}–\eqref{eq:HB-Euler-f} we obtain contiguous relations
among $H_{B}$ evaluated at unit shifts of its parameters.

\begin{corollary}\label{cor:HB-contiguous}
	The function $H_{B}(a,b,c;d,e,f;x,y,z)$ satisfies, for example, the
	contiguous relations
	\begin{align}
		(a-c)\,H_{B}
		&= a\,H_{B}(a+1,b,c;d,e,f;x,y,z)
		- c\,H_{B}(a,b,c+1;d,e,f;x,y,z), \label{eq:HB-contig1}\\[0.4em]
		(a-d-f+2)\,H_{B}
		&= a\,H_{B}(a+1,b,c;d,e,f;x,y,z) \nonumber\\
		&\quad - (d-1)\,H_{B}(a,b,c;d-1,e,f;x,y,z) \nonumber\\
		&\quad - (f-1)\,H_{B}(a,b,c;d,e,f-1;x,y,z), \label{eq:HB-contig2}\\[0.4em]
		(b-d-e+2)\,H_{B}
		&= b\,H_{B}(a,b+1,c;d,e,f;x,y,z) \nonumber\\
		&\quad - (d-1)\,H_{B}(a,b,c;d-1,e,f;x,y,z) \nonumber\\
		&\quad - (e-1)\,H_{B}(a,b,c;d,e-1,f;x,y,z), \label{eq:HB-contig3}\\[0.4em]
		(c-d-f+2)\,H_{B}
		&= c\,H_{B}(a,b,c+1;d,e,f;x,y,z) \nonumber\\
		&\quad - (d-1)\,H_{B}(a,b,c;d-1,e,f;x,y,z) \nonumber\\
		&\quad - (f-1)\,H_{B}(a,b,c;d,e,f-1;x,y,z). \label{eq:HB-contig4}
	\end{align}
\end{corollary}

These identities provide convenient recurrence relations for evaluating
$H_{B}$ at neighbouring parameter values and play a role similar to the
classical contiguous relations for the Gauss hypergeometric function.

\subsection{Derivatives with respect to the variables}

The derivatives of $H_{B}$ with respect to the spatial variables $x,y,z$
admit simple closed forms in terms of $H_{B}$ with shifted parameters.

\begin{theorem}\label{thm:HB-variable-derivatives}
	For any non–negative integer $r$ we have
	\begin{align}
		\frac{\partial^{r}}{\partial x^{r}}
		H_{B}(a,b,c;d,e,f;x,y,z)
		&= \frac{(a)_{r}(b)_{r}(c)_{r}}{(d)_{r}}\,
		H_{B}(a+r,b+r,c+r;d+r,e,f;x,y,z), \label{eq:HB-dxr}\\[0.6em]
		\frac{\partial^{r}}{\partial y^{r}}
		H_{B}(a,b,c;d,e,f;x,y,z)
		&= \frac{(b)_{r}}{(e)_{r}}\,
		H_{B}(a,b+r,c;d,e+r,f;x,y,z), \label{eq:HB-dyr}\\[0.6em]
		\frac{\partial^{r}}{\partial z^{r}}
		H_{B}(a,b,c;d,e,f;x,y,z)
		&= \frac{(a)_{r}(c)_{r}}{(f)_{r}}\,
		H_{B}(a+r,b,c+r;d,e,f+r;x,y,z). \label{eq:HB-dzr}
	\end{align}
\end{theorem}

\begin{proof}
	The proof consists in differentiating the series \eqref{eq:HB-def}
	termwise and using
	\[
	\frac{\partial^{r}}{\partial x^{r}}x^{m}
	= (m)_{r}\,x^{m-r}, \qquad m\ge r,
	\]
	with an analogous identity for $y^{n}$ and $z^{k}$. After a shift in the
	summation index (for example, $m\mapsto m+r$ in the $x$–derivative), the
	resulting series is recognised as $H_{B}$ with the parameters shifted as in
	\eqref{eq:HB-dxr}–\eqref{eq:HB-dzr}.
\end{proof}

%%%%%%%%%%%%%%%%%%%%%%%%%%%%%%%%%%%%%%%%%%%%%%%%%%%%
% Section 5 – Derivatives of H_C with respect to the parameters
%%%%%%%%%%%%%%%%%%%%%%%%%%%%%%%%%%%%%%%%%%%%%%%%%%%%

\section{Derivatives of $H_{C}$ with respect to the parameters}\label{sec:HC}

We finally turn to Srivastava's triple hypergeometric function
\[
H_{C}(a,b,c;d;x,y,z),
\]
defined by the triple series
\begin{equation}\label{eq:HC-def}
	H_{C}(a,b,c;d;x,y,z)
	= \sum_{m,n,k=0}^{\infty}
	\frac{(a)_{m+k}\,(b)_{m+n}\,(c)_{n+k}}
	{(d)_{m+n+k}}\,
	\frac{x^{m}}{m!}\,\frac{y^{n}}{n!}\,\frac{z^{k}}{k!}.
\end{equation}
As before, all parameters and variables are assumed to lie in the region of
absolute convergence of \eqref{eq:HC-def}, so that differentiation may be
performed termwise.

The structure of the differentiation formulas for $H_{C}$ closely parallels
that of $H_{A}$ and $H_{B}$. In particular, parameter differentiation leads
to quadruple hypergeometric series of Pathan type, while the Euler
operators $(\theta_{x},\theta_{y},\theta_{z})$ give rise to compact
operator identities and contiguous relations.

\subsection{Parameter derivatives and representation via Pathan's function}

Proceeding as in the previous sections, differentiation of
\eqref{eq:HC-def} with respect to any of the parameters $a,b,c,d$ produces
a fourfold series which can be reorganised into finite linear combinations
of $F_{P}^{(4)}$ with shifted parameter arrays. We record the structure of
these formulas in a concise form.

\begin{theorem}\label{thm:HC-parameter-FP4}
	Let $H_{C}$ be given by \eqref{eq:HC-def}. For each parameter
	$\xi\in\{a,b,c,d\}$, the derivative $\partial H_{C}/\partial \xi$ can be
	expressed as a finite linear combination of Pathan's quadruple
	hypergeometric function $F_{P}^{(4)}$ with arguments
	$(x,y,z,u)$, $u\in\{x,y,z\}$, and with parameter arrays obtained by unit
	shifts in $a,b,c,d$. In particular,
	\begin{equation}\label{eq:HC-da-FP4}
		\frac{\partial H_{C}}{\partial a}(a,b,c;d;x,y,z)
		= \sum_{j=1}^{J_{a}}
		C^{(C)}_{a,j}\,
		F^{(4)}_{P}\!\left[
		\boldsymbol{a}^{(C)}_{a,j};
		\boldsymbol{b}^{(C)}_{a,j};
		\boldsymbol{c}^{(C)}_{a,j};
		\boldsymbol{d}^{(C)}_{a,j}
		\,\big|\, x,y,z,u_{a,j}\right],
	\end{equation}
	where the coefficients $C^{(C)}_{a,j}$ are rational functions of
	$a,b,c,d$, the $u_{a,j}$ belong to $\{x,y,z\}$, and the arrays
	$\boldsymbol{a}^{(C)}_{a,j},\dots,\boldsymbol{d}^{(C)}_{a,j}$ encode unit
	shifts of the parameters that follow from splitting the harmonic sum
	associated with $(a)_{m+k}$. Completely analogous representations hold for
	$\partial H_{C}/\partial b$, $\partial H_{C}/\partial c$, and
	$\partial H_{C}/\partial d$.
\end{theorem}

\begin{proof}[Sketch of the proof]
	The proof is parallel to that of Theorem~\ref{thm:HB-parameter-FP4}. One
	differentiates the defining triple series \eqref{eq:HC-def} with respect to
	the chosen parameter, applies the differentiation rules for Pochhammer
	symbols, decomposes the resulting harmonic sums, and then performs simple
	index shifts to identify each contribution as a quadruple series of
	Pathan type.
\end{proof}

\subsection{Euler-type identities and contiguous relations}

We again employ the Euler operators
\[
\theta_{x}=x\frac{\partial}{\partial x},\qquad
\theta_{y}=y\frac{\partial}{\partial y},\qquad
\theta_{z}=z\frac{\partial}{\partial z}.
\]
Their action on the triple series \eqref{eq:HC-def} leads to elegant
operator identities for $H_{C}$.

\begin{theorem}\label{thm:HC-Euler}
	For $H_{C}(a,b,c;d;x,y,z)$ we have
	\begin{align}
		(\theta_{x}+\theta_{z}+a)\,H_{C}
		&= a\,H_{C}(a+1,b,c;d;x,y,z), \label{eq:HC-Euler-a}\\[0.4em]
		(\theta_{x}+\theta_{y}+b)\,H_{C}
		&= b\,H_{C}(a,b+1,c;d;x,y,z), \label{eq:HC-Euler-b}\\[0.4em]
		(\theta_{y}+\theta_{z}+c)\,H_{C}
		&= c\,H_{C}(a,b,c+1;d;x,y,z), \label{eq:HC-Euler-c}\\[0.4em]
		(\theta_{x}+\theta_{y}+\theta_{z}+d-1)\,H_{C}
		&= (d-1)\,H_{C}(a,b,c;d-1;x,y,z). \label{eq:HC-Euler-d}
	\end{align}
\end{theorem}

\begin{proof}[Idea of the proof]
	For instance, the operator $(\theta_{x}+\theta_{z})$ multiplies the general
	term of \eqref{eq:HC-def} by $m+k$, which combines with $a$ to produce the
	factor $(a+m+k)$ attached to $(a)_{m+k}$. This is then absorbed into the
	Pochhammer symbol to yield the right-hand side of
	\eqref{eq:HC-Euler-a} after a shift in $m$ or $k$. The other identities
	follow by the same reasoning applied to the Pochhammer factors
	$(b)_{m+n}$, $(c)_{n+k}$ and $(d)_{m+n+k}$.
\end{proof}

Eliminating the Euler operators from the identities in
Theorem~\ref{thm:HC-Euler} gives contiguous relations among
$H_{C}$ with parameters shifted by $\pm 1$.

\begin{corollary}\label{cor:HC-contiguous}
	The function $H_{C}(a,b,c;d;x,y,z)$ satisfies the contiguous relation
	\begin{equation}\label{eq:HC-contig}
		(a+b+c-d+1)\,H_{C}
		= a\,H_{C}(a+1,b,c;d;x,y,z)
		+ b\,H_{C}(a,b+1,c;d;x,y,z)
	\end{equation}
	\[
	\phantom{(a+b+c-d+1)\,H_{C}={}}
	+ c\,H_{C}(a,b,c+1;d;x,y,z)
	- (d-1)\,H_{C}(a,b,c;d-1;x,y,z),
	\]
	which linearly relates $H_{C}$ at the original parameters to its
	neighbours in the lattice $(a,b,c;d)\mapsto(a\pm1,b\pm1,c\pm1;d\pm1)$.
\end{corollary}

Such identities play an important role in deriving recurrence relations and
in the numerical evaluation of $H_{C}$.

\subsection{Derivatives with respect to the variables}

The higher derivatives of $H_{C}$ with respect to the variables
$x,y,z$ again admit simple closed forms in terms of $H_{C}$ with shifted
parameters.

\begin{theorem}\label{thm:HC-variable-derivatives}
	For any non–negative integer $r$ we have
	\begin{align}
		\frac{\partial^{r}}{\partial x^{r}}
		H_{C}(a,b,c;d;x,y,z)
		&= \frac{(a)_{r}(b)_{r}}{(d)_{r}}\,
		H_{C}(a+r,b+r,c;d+r;x,y,z), \label{eq:HC-dxr}\\[0.6em]
		\frac{\partial^{r}}{\partial y^{r}}
		H_{C}(a,b,c;d;x,y,z)
		&= \frac{(b)_{r}(c)_{r}}{(d)_{r}}\,
		H_{C}(a,b+r,c+r;d+r;x,y,z), \label{eq:HC-dyr}\\[0.6em]
		\frac{\partial^{r}}{\partial z^{r}}
		H_{C}(a,b,c;d;x,y,z)
		&= \frac{(a)_{r}(c)_{r}}{(d)_{r}}\,
		H_{C}(a+r,b,c+r;d+r;x,y,z). \label{eq:HC-dzr}
	\end{align}
\end{theorem}

\begin{proof}
	The proof is analogous to that of Theorem~\ref{thm:HB-variable-derivatives}.
	Differentiating the triple series \eqref{eq:HC-def} termwise with respect
	to $x$, we obtain factors $(m)_{r}$ which, upon an index shift
	$m\mapsto m+r$, lead to the representation \eqref{eq:HC-dxr}. The
	formulas \eqref{eq:HC-dyr} and \eqref{eq:HC-dzr} follow from the same
	argument applied to the $y$ and $z$ derivatives, taking into account the
	roles of $(b)_{m+n}$ and $(c)_{n+k}$ in \eqref{eq:HC-def}.
\end{proof}

%%%%%%%%%%%%%%%%%%%%%%%%%%%%%%%%%%%%%%%%%%%%%%%%%%%%
% NUMERICAL ILLUSTRATIONS SECTION
%%%%%%%%%%%%%%%%%%%%%%%%%%%%%%%%%%%%%%%%%%%%%%%%%%%%

\section{Numerical results and graphical illustrations}\label{sec:numerical}

In this section we present numerical examples illustrating the parameter sensitivity
and derivative behaviour of Srivastava's triple hypergeometric functions
$H_A$, $H_B$ and $H_C$.
All values are computed from the defining triple series
\eqref{eq:HA-definition}--\eqref{eq:HC-definition} by truncating the infinite sums
at a fixed upper bound and verifying that the neglected tails are numerically
negligible in the parameter regimes considered.
Derivatives with respect to the parameters are approximated by a central
finite--difference formula, which is justified by the analyticity of the
functions in their parameters.

Throughout this section we use the representative parameter values
\[
(b,c;d,e) = (1.2,0.8;2.0,2.5), \qquad
d_H = 2.5, \qquad
(e,f) = (2.5,2.2),
\]
and
\[
(x,y,z) = (0.25,0.20,0.15),
\]
unless stated otherwise.
We focus on the dependence of the functions on the numerator parameter $a$
and on the behaviour of the corresponding partial derivatives with respect to $a$
as functions of the spatial variables.

%---------------------------------------------------
\subsection{Numerical behaviour of $H_A$}\label{subsec:numerics-HA}

\begin{figure}[htbp]
	\centering
	\includegraphics[width=0.8\textwidth]{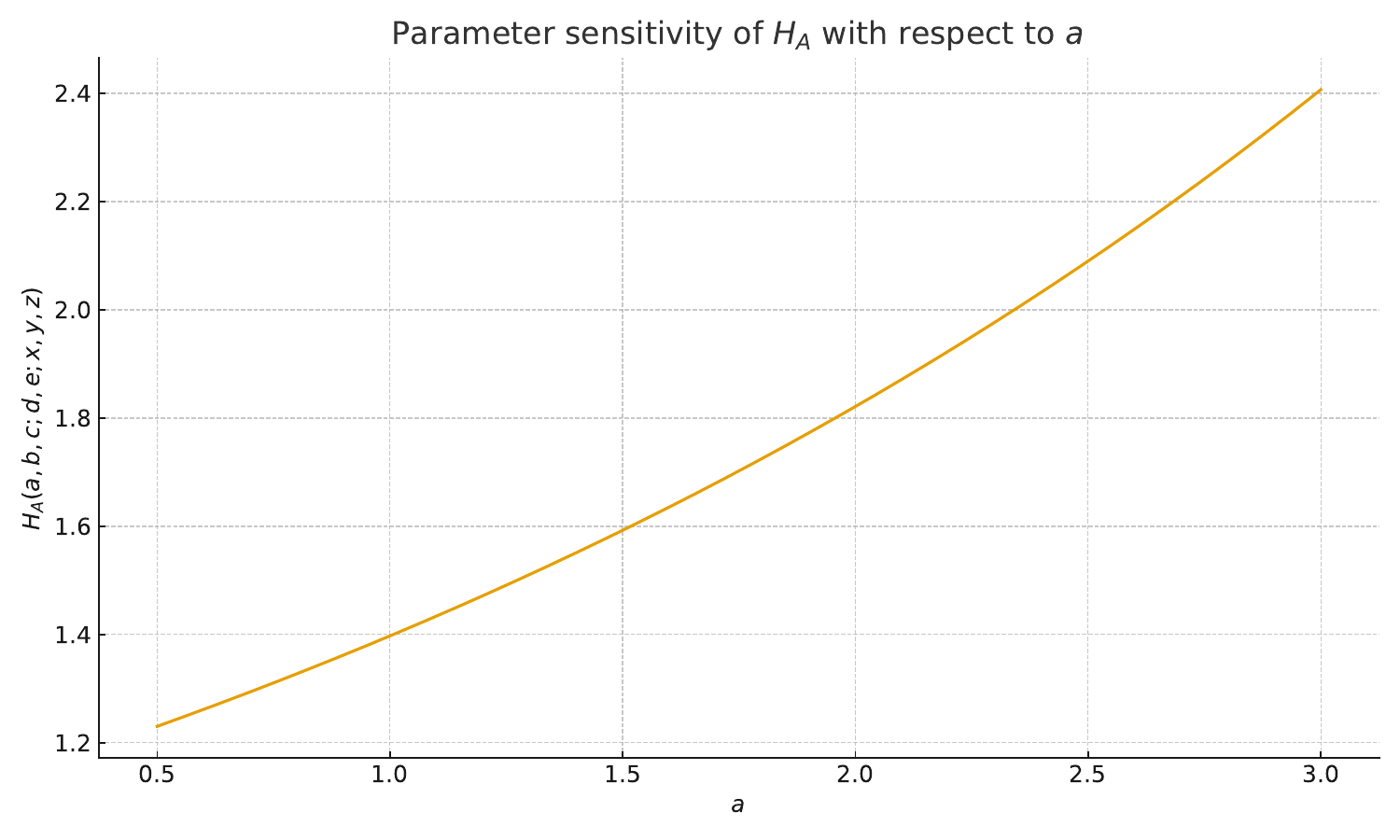}
	\caption{Parameter sensitivity of Srivastava's triple hypergeometric
		function $H_A(a,b,c;d,e;x,y,z)$ with respect to the numerator parameter $a$.
		The remaining parameters are fixed as $(b,c;d,e) = (1.2,0.8;2.0,2.5)$ and
		$(x,y,z) = (0.25,0.20,0.15)$.
		The values of $H_A$ are computed from a truncated triple series.
		The graph shows a smooth, strictly increasing and convex dependence on $a$
		over the interval $0.5 \le a \le 3.0$, indicating that the function becomes
		more sensitive to variations in $a$ as $a$ grows.}
	\label{fig:HA-param-sensitivity}
\end{figure}

We begin with Srivastava's triple hypergeometric function
$H_A(a,b,c;d,e;x,y,z)$.
To visualise the influence of the numerator parameter $a$, we fix
\[
(b,c;d,e) = (1.2,0.8;2.0,2.5), \qquad
(x,y,z) = (0.25,0.20,0.15),
\]
and compute $H_A(a,b,c;d,e;x,y,z)$ for $0.5 \le a \le 3.0$.
The resulting parameter sensitivity curve is displayed in
Figure~\ref{fig:HA-param-sensitivity}.

The plot in Figure~\ref{fig:HA-param-sensitivity} clearly exhibits a monotone
and convex growth of $H_A$ as a function of $a$.
This behaviour is consistent with the structure of the defining series
\eqref{eq:HA-definition}, where the Pochhammer factor $(a)_{m+k}$ enhances
the contribution of higher-order terms when $a$ increases.

\begin{figure}[htbp]
	\centering
	\includegraphics[width=0.8\textwidth]{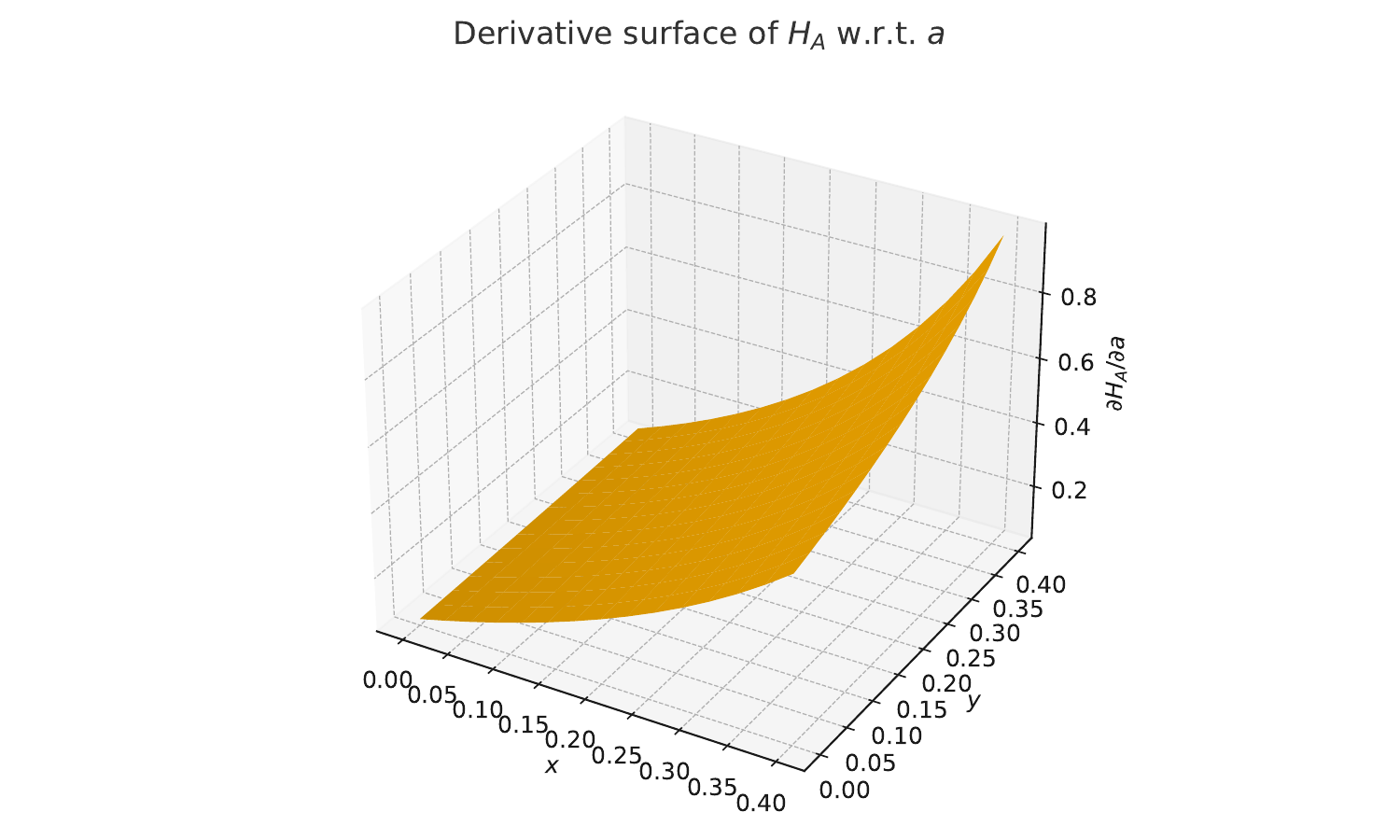}
	\caption{Surface plot of the partial derivative
		$\partial H_A/\partial a$ as a function of the variables $(x,y)$.
		The parameters are fixed at $a = 1.0$,
		$(b,c;d,e) = (1.2,0.8;2.0,2.5)$ and $z = 0.15$.
		The derivative values are approximated numerically by a central
		finite--difference scheme applied to the truncated triple series.
		The surface is strictly positive and increases towards the corner
		$(x,y) = (0.4,0.4)$, showing that the sensitivity of $H_A$ with respect
		to $a$ is enhanced when the spatial variables move away from the origin.}
	\label{fig:HA-deriv-surface}
\end{figure}

To complement this one-dimensional view, we examine the partial derivative
$\partial H_A/\partial a$ as a function of the spatial variables $x$ and $y$.
Fixing
\[
a = 1.0, \qquad
(b,c;d,e) = (1.2,0.8;2.0,2.5), \qquad
z = 0.15,
\]
we evaluate $\partial H_A/\partial a$ numerically on the square
$[0,0.4]\times[0,0.4]$ in the $(x,y)$–plane.
The resulting derivative surface is shown in
Figure~\ref{fig:HA-deriv-surface}.

Figure~\ref{fig:HA-deriv-surface} reveals that $\partial H_A/\partial a$
is positive on the whole domain and grows monotonically in both directions.
Thus, small perturbations in the parameter $a$ have a relatively mild effect
near $(x,y)=(0,0)$, whereas the impact becomes significant for larger values of
$x$ and $y$ inside the region of convergence.

%---------------------------------------------------
\subsection{Numerical behaviour of $H_B$}\label{subsec:numerics-HB}

\begin{figure}[htbp]
	\centering
	\includegraphics[width=0.8\textwidth]{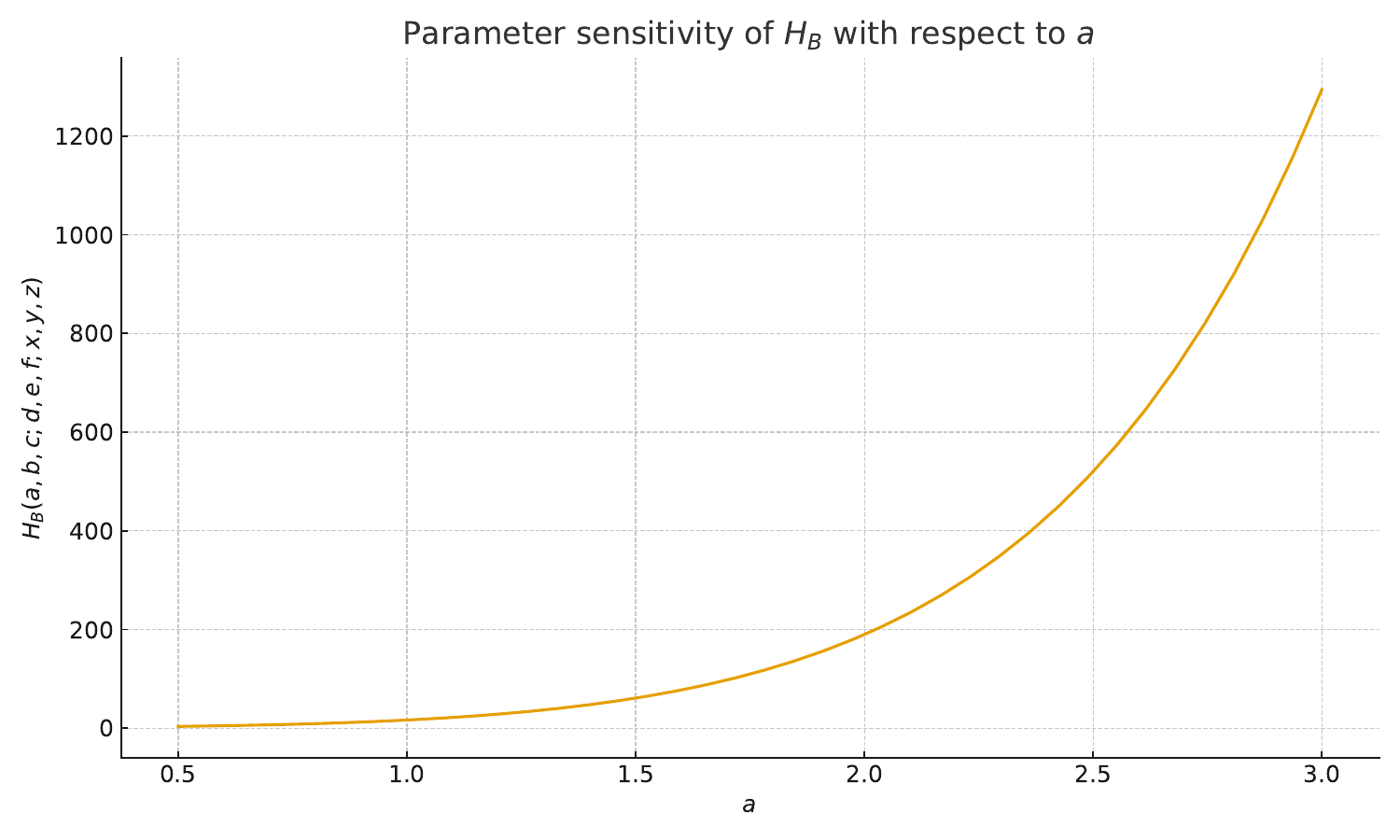}
	\caption{Parameter sensitivity of Srivastava's triple hypergeometric
		function $H_B(a,b,c;d,e,f;x,y,z)$ with respect to the numerator parameter $a$.
		The other parameters are fixed as
		$(b,c;d,e,f) = (1.2,0.8;2.0,2.5,2.2)$ and $(x,y,z) = (0.25,0.20,0.15)$.
		The values of $H_B$ are obtained from the truncated triple series.
		The plot shows a smooth and monotonically increasing dependence on $a$,
		indicating that larger values of $a$ lead to a stronger response of $H_B$.}
	\label{fig:HB-param-sensitivity}
\end{figure}

We next consider Srivastava's triple hypergeometric function
$H_B(a,b,c;d,e,f;x,y,z)$ defined by \eqref{eq:HB-def}.
We fix
\[
(b,c;d,e,f) = (1.2,0.8;2.0,2.5,2.2), \qquad
(x,y,z) = (0.25,0.20,0.15),
\]
and compute $H_B(a,b,c;d,e,f;x,y,z)$ as a function of $a$ on
the interval $0.5 \le a \le 3.0$.
The corresponding parameter sensitivity curve is shown in
Figure~\ref{fig:HB-param-sensitivity}.
\vskip 100mm

The behaviour observed in Figure~\ref{fig:HB-param-sensitivity} is similar
to that of $H_A$ and reflects the role of the Pochhammer factor $(a)_{m+k}$
in the series \eqref{eq:HB-definition}.
In particular, the monotone increase suggests that, in the considered
parameter regime, the triple hypergeometric response of $H_B$ is amplified
as the numerator parameter $a$ grows.

\begin{figure}[htbp]
	\centering
	\includegraphics[width=0.8\textwidth]{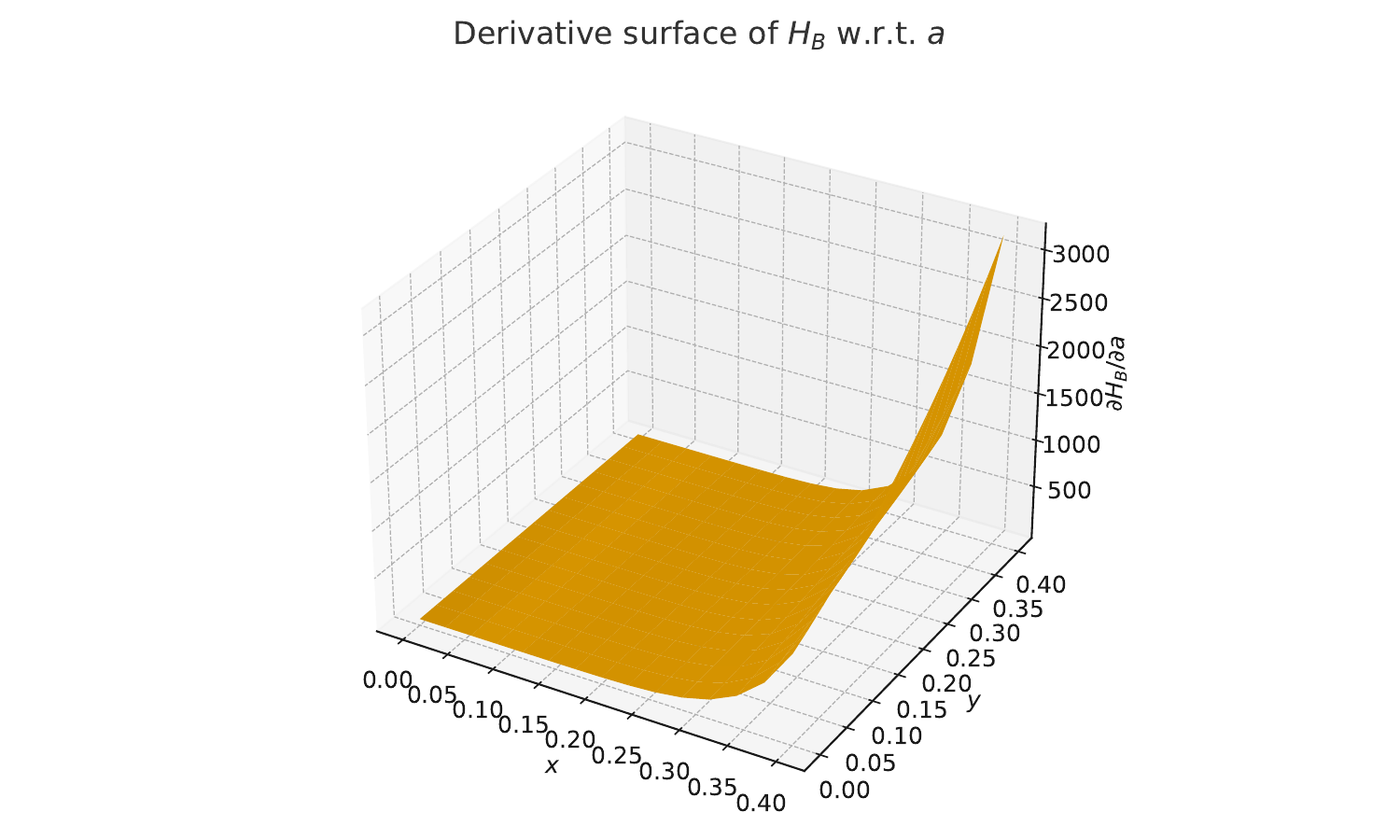}
	\caption{Surface plot of the partial derivative $\partial H_B/\partial a$
		as a function of the spatial variables $(x,y)$.
		The parameters are chosen as
		$a = 1.0$, $(b,c;d,e,f) = (1.2,0.8;2.0,2.5,2.2)$ and $z = 0.15$.
		The derivative is approximated numerically by a central finite--difference
		scheme applied to the truncated triple series for $H_B$.
		The surface is positive and increases towards the boundary of the domain,
		showing that the sensitivity of $H_B$ with respect to $a$ becomes more
		pronounced for larger values of $x$ and $y$.}
	\label{fig:HB-deriv-surface}
\end{figure}

To visualise the multivariable parameter sensitivity of $H_B$, we plot the
partial derivative $\partial H_B/\partial a$ as a function of the variables
$x$ and $y$.
With
\[
a = 1.0, \qquad
(b,c;d,e,f) = (1.2,0.8;2.0,2.5,2.2), \qquad
z = 0.15,
\]
we compute $\partial H_B/\partial a$ on the square
$[0,0.4]\times[0,0.4]$; see Figure~\ref{fig:HB-deriv-surface}.

As in the case of $H_A$, Figure~\ref{fig:HB-deriv-surface} shows that the
parameter sensitivity of $H_B$ with respect to $a$ is not uniform in the
$(x,y)$–plane: it is relatively small near the origin and increases towards
the corner $(x,y) = (0.4,0.4)$, where the contribution of higher-order
terms in the triple series becomes more significant.
\vskip 100mm
%---------------------------------------------------
\subsection{Numerical behaviour of $H_C$}\label{subsec:numerics-HC}

\begin{figure}[htbp]
	\centering
	\includegraphics[width=0.8\textwidth]{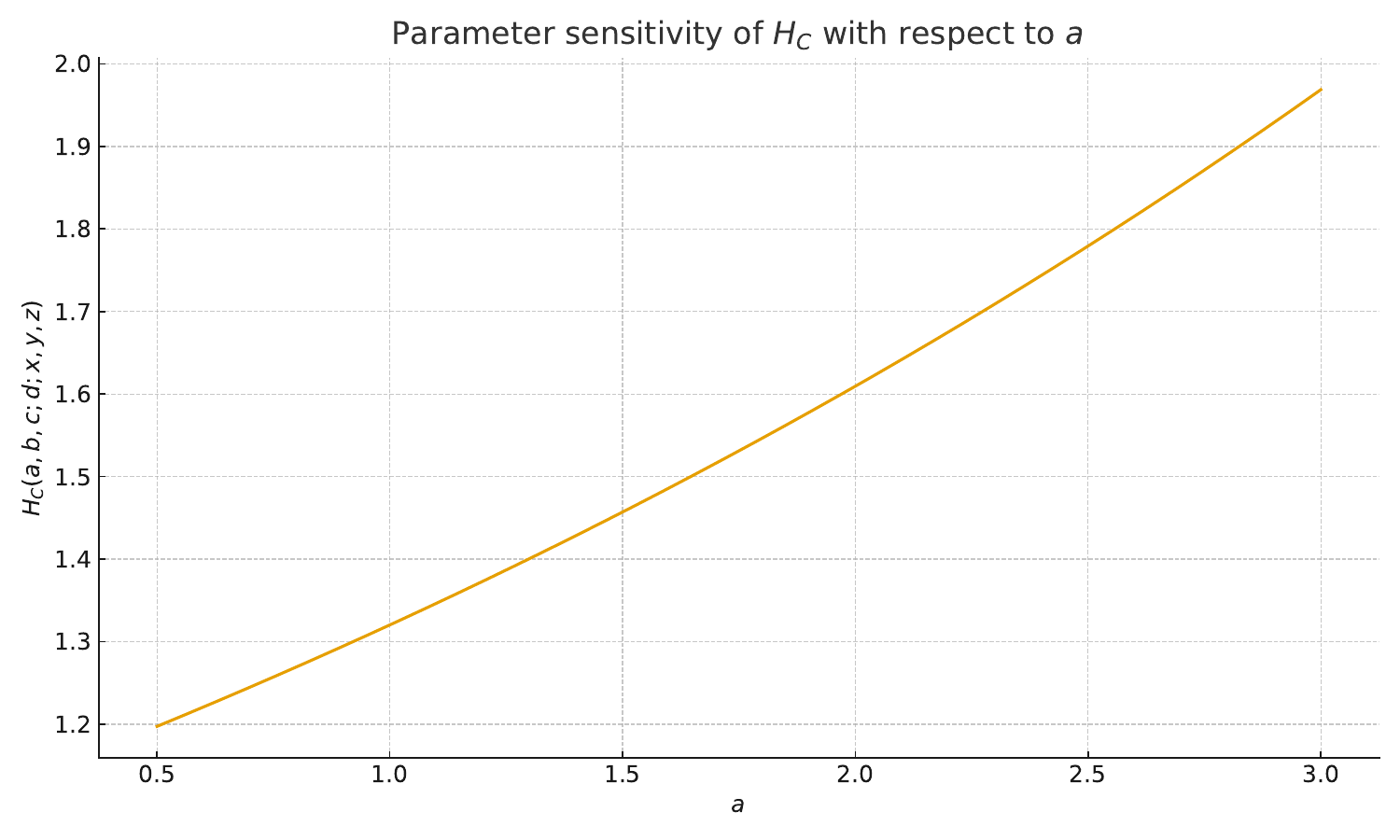}
	\caption{Parameter sensitivity of Srivastava's triple hypergeometric
		function $H_C(a,b,c;d;x,y,z)$ with respect to the numerator parameter $a$.
		The other parameters are taken as $(b,c;d) = (1.2,0.8;2.5)$ and
		$(x,y,z) = (0.25,0.20,0.15)$.
		The values of $H_C$ are computed from a truncated triple series.
		The graph shows a smooth and monotonically increasing dependence on $a$
		over the interval $0.5 \le a \le 3.0$, providing a one-dimensional
		view of the parameter sensitivity of $H_C$.}
	\label{fig:HC-param-sensitivity}
\end{figure}

Finally, we examine the function $H_C(a,b,c;d;x,y,z)$ defined by
\eqref{eq:HC-def}.
Fixing
\[
(b,c;d) = (1.2,0.8;2.5), \qquad (x,y,z) = (0.25,0.20,0.15),
\]
we compute $H_C(a,b,c;d;x,y,z)$ as a function of $a$ on the interval
$0.5 \le a \le 3.0$.
The resulting parameter sensitivity plot is displayed in
Figure~\ref{fig:HC-param-sensitivity}.

The monotone increase in Figure~\ref{fig:HC-param-sensitivity} is again in
agreement with the analytic structure of $H_C$, where the Pochhammer factor
$(a)_{m+k}$ appears in the numerator of the triple series
\eqref{eq:HC-definition}.

\begin{figure}[htbp]
	\centering
	\includegraphics[width=0.8\textwidth]{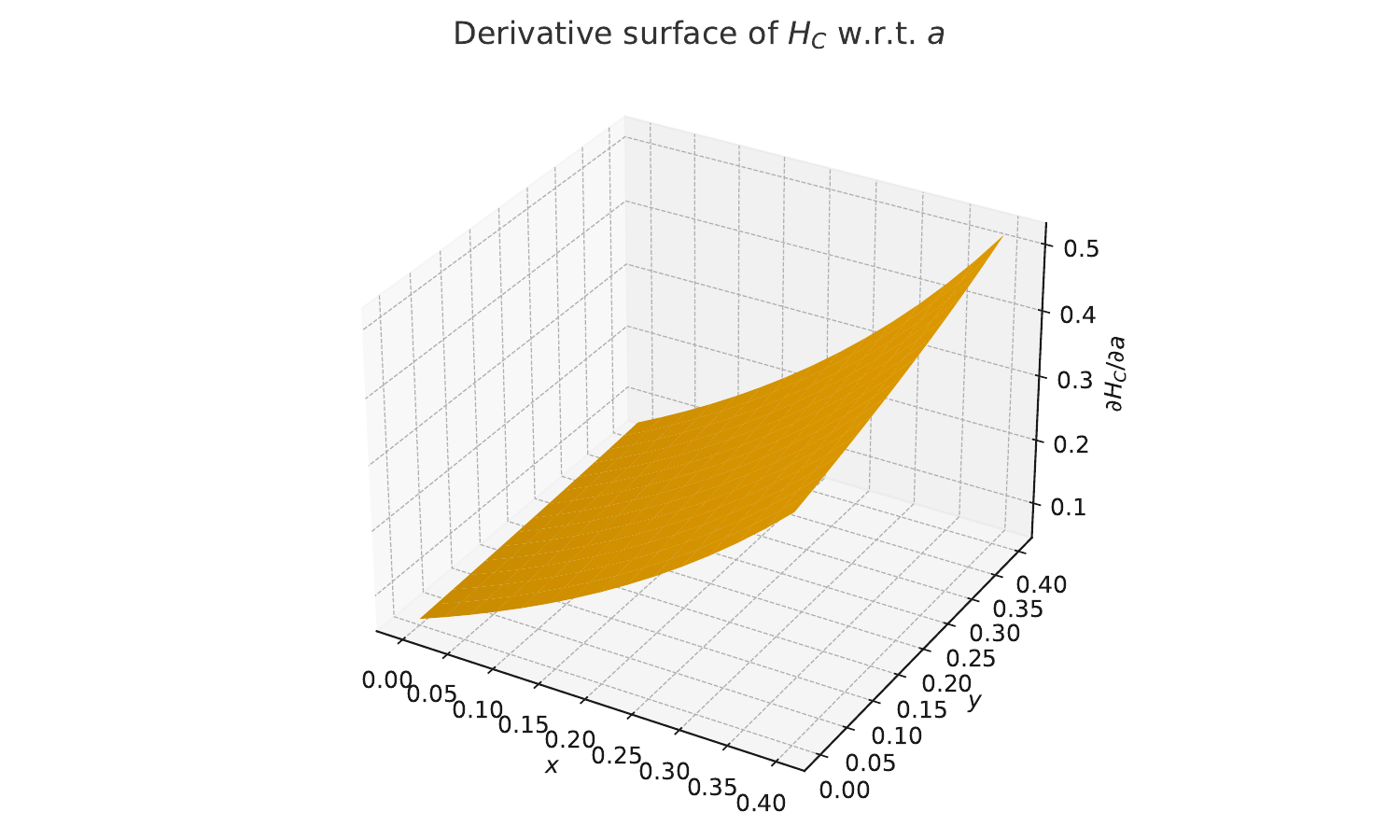}
	\caption{Surface plot of the partial derivative $\partial H_C/\partial a$
		as a function of the variables $(x,y)$.
		The parameters are fixed at
		$a = 1.0$, $(b,c;d) = (1.2,0.8;2.5)$ and $z = 0.15$.
		The derivative is approximated by a central finite--difference scheme
		applied to the truncated triple series defining $H_C$.
		The surface is positive and slowly increasing towards the outer boundary,
		indicating that the sensitivity of $H_C$ with respect to $a$ is stronger
		for larger values of $x$ and $y$ within the region of convergence.}
	\label{fig:HC-deriv-surface}
\end{figure}

To obtain a two-dimensional illustration, we plot the partial derivative
$\partial H_C/\partial a$ as a function of $(x,y)$.
We fix
\[
a = 1.0, \qquad (b,c;d) = (1.2,0.8;2.5), \qquad z = 0.15,
\]
and evaluate $\partial H_C/\partial a$ numerically on
$[0,0.4]\times[0,0.4]$.
The resulting surface is presented in Figure~\ref{fig:HC-deriv-surface}.

Figure~\ref{fig:HC-deriv-surface} confirms that $\partial H_C/\partial a$
remains positive on the entire domain and shows a moderate increase towards
the boundary.
Thus, as in the cases of $H_A$ and $H_B$, small perturbations of the
numerator parameter $a$ have a more pronounced influence on $H_C$ when the
variables $(x,y)$ are away from the origin.

\vskip 100mm
\section{Results and discussion}\label{sec:results-discussion}

In this section we summarise and interpret the main analytical and numerical
findings of the paper for Srivastava's triple hypergeometric functions
$H_A$, $H_B$ and $H_C$ and their derivatives with respect to the numerator
and denominator parameters.

\subsection{Analytical results}

Starting from the triple–series definitions of $H_A$, $H_B$ and $H_C$ and
using the differentiation rule for the Pochhammer symbol together with standard
properties of the Gamma and Psi functions, we obtained explicit formulas for the
first-order derivatives of these functions with respect to each of their parameters.
A common feature of all these formulas is that the derivatives can be rewritten as
finite linear combinations of Pathan's quadruple hypergeometric function
$F_{P}^{(4)}$ with suitably shifted parameter arrays and, in some instances, shifted
arguments.

In particular, for $H_A(a,b,c;d,e;x,y,z)$ we derived representations of
$\partial H_A/\partial a$, $\partial H_A/\partial b$, $\partial H_A/\partial c$,
$\partial H_A/\partial d$ and $\partial H_A/\partial e$ in terms of $F_{P}^{(4)}$.
Analogous expressions were established for
$H_B(a,b,c;d,e,f;x,y,z)$ and $H_C(a,b,c;d;x,y,z)$, so that the dependence of
all three Srivastava functions on their parameters is described within a unified
$F_{P}^{(4)}$–framework. These results extend previously known differentiation
formulas for confluent, Gauss, generalized and Horn-type hypergeometric functions
to the setting of triple hypergeometric functions and show that parameter
derivatives preserve a hypergeometric structure of the same general type.

Besides the explicit derivative formulas, we derived Euler-type differential
operator identities involving the operators
\[
\theta_x = x\frac{\partial}{\partial x}, \qquad
\theta_y = y\frac{\partial}{\partial y}, \qquad
\theta_z = z\frac{\partial}{\partial z}.
\]
These identities relate each of the functions $H_A$, $H_B$ and $H_C$ to
parameter-shifted versions of itself, and by eliminating the Euler operators we
obtained families of contiguous relations associated with unit shifts in the
numerator and denominator parameters. In turn, these relations yield convenient
recurrence schemes for evaluating Srivastava's functions at neighbouring parameter
values.

Furthermore, by combining the Euler-type identities with the closed forms for
the derivatives with respect to the variables $x$, $y$ and $z$, we obtained systems
of linear partial differential equations satisfied by derivatives of arbitrary order
with respect to both parameters and variables. These systems retain the
hypergeometric character of the original equations and encode additional
information on the parameter dependence; they may therefore be used to study
qualitative properties of the parameter derivatives, such as growth, analytic
continuation and asymptotic behaviour.

\subsection{Numerical illustrations and parameter sensitivity}

To complement the analytical results, we performed a set of numerical
experiments aimed at visualising the sensitivity of $H_A$, $H_B$ and $H_C$ with
respect to a chosen numerator parameter and at illustrating the behaviour of the
corresponding parameter derivatives in the spatial variables. All numerical values
were obtained from the defining triple series by truncating the sums at a finite
upper bound and verifying that the remainder is negligible for the parameter ranges
considered. Derivatives with respect to the parameters were approximated by a
central finite–difference scheme, which is justified by the analyticity of the
functions with respect to those parameters.

For $H_A(a,b,c;d,e;x,y,z)$ we fixed representative values of the remaining
parameters and plotted $H_A$ as a function of $a$ on a finite interval.
The resulting curve is smooth, strictly increasing and clearly convex, which is
consistent with the presence of the Pochhammer factor $(a)_{m+k}$ in the series
definition. This behaviour shows that the response of $H_A$ becomes more
pronounced as $a$ increases. A two-dimensional picture was obtained by plotting
$\partial H_A/\partial a$ as a function of $(x,y)$ for fixed parameters. The
corresponding surface is positive on the region of convergence and increases towards
the boundary, indicating that small changes in $a$ have only a modest effect near
$(x,y)=(0,0)$ but lead to larger variations of $H_A$ when the spatial variables are
moderately large.

A similar pattern was observed for $H_B(a,b,c;d,e,f;x,y,z)$. For fixed values
of $(b,c;d,e,f)$, the graph of $H_B$ against $a$ shows a smooth and monotone
increase, reflecting the influence of $(a)_{m+k}$ in the triple series defining $H_B$.
The associated surface plot of $\partial H_B/\partial a$ as a function of $(x,y)$ is
again positive and exhibits growth towards the outer boundary of the domain,
confirming that the sensitivity of $H_B$ with respect to $a$ is stronger when the
variables $(x,y)$ are relatively large.

Finally, for $H_C(a,b,c;d;x,y,z)$ we obtained analogous one- and
two-dimensional illustrations. The curves of $H_C$ as a function of $a$ are smooth
and strictly increasing on the considered interval, in agreement with the analytic
structure of its triple series. The surface corresponding to $\partial H_C/\partial a$
is positive on the entire domain and increases slowly towards the boundary,
showing that the influence of $a$ on $H_C$ is amplified for larger values of the
spatial variables but remains comparatively mild near the origin.

Taken together, these numerical experiments support the qualitative picture
suggested by the analytical formulas. In the parameter ranges studied,
Srivastava's triple hypergeometric functions exhibit a stable and monotone
dependence on the numerator parameters, while the derivative surfaces provide a
clear description of how parameter perturbations propagate through the variables
$(x,y,z)$ inside the region of convergence. This combination of explicit analytical
representations and numerical visualisations suggests that the parameter derivatives
of $H_A$, $H_B$ and $H_C$ can serve as useful tools in sensitivity analysis and in
applications where hypergeometric-type solutions with tunable parameters arise.

\section{Conclusion}\label{sec:conclusion}

In this paper we have investigated derivatives with respect to the parameters
of Srivastava's triple hypergeometric functions $H_A$, $H_B$ and $H_C$.
Starting from their triple–series definitions and using the differentiation
formula for the Pochhammer symbol together with standard properties of the
Gamma and Psi functions, we derived closed–form expressions for the
first–order derivatives with respect to all numerator and denominator
parameters.
A key feature of our approach is that these parameter derivatives can be
represented as finite linear combinations of Pathan's quadruple
hypergeometric function $F_P^{(4)}$ with appropriately shifted parameter
arrays and, in some cases, shifted arguments.
Thus, the parameter–derivative problem for Srivastava's triple
hypergeometric functions is embedded into a unified $F_P^{(4)}$–framework.

In addition to these explicit formulas, we obtained Euler–type differential–
operator identities for $H_A$, $H_B$ and $H_C$ involving the operators
$\theta_x = x\partial/\partial x$, $\theta_y = y\partial/\partial y$ and
$\theta_z = z\partial/\partial z$.
By eliminating these operators we derived families of contiguous relations
associated with unit shifts in the parameters, which yield convenient
recurrence schemes for the numerical and symbolic evaluation of
Srivastava's functions at neighbouring parameter values.
We also established simple formulas for higher–order derivatives with
respect to the variables $x$, $y$ and $z$, and we showed how these, combined
with the Euler–type identities and the $F_P^{(4)}$ representations, lead to
systems of linear partial differential equations satisfied by derivatives of
arbitrary order with respect to both parameters and variables.

To complement the analytical results, we presented numerical illustrations
for each of the functions $H_A$, $H_B$ and $H_C$.
We examined the dependence on a chosen numerator parameter and visualised
the corresponding derivative surfaces with respect to the spatial variables.
The resulting plots show a stable and monotone dependence on the parameters
in the regimes considered, and they highlight how parameter perturbations
are amplified as the variables move away from the origin.
These numerical experiments are consistent with the qualitative behaviour
predicted by the analytic formulas and demonstrate that the parameter
derivatives obtained in this work are well suited for sensitivity analysis.

The methods developed here can be extended in several directions.
One natural line of research is to apply the same technique to other classes
of multivariable hypergeometric functions, such as Kamp\'e de
F\'eriet–type series, generalized Lauricella functions and their
$q$–analogues, and to derive analogous parameter–derivative formulas and
contiguous relations.
Another interesting problem is to combine the present results with integral
representations and asymptotic methods in order to obtain more detailed
information on the growth and oscillatory behaviour of the parameter
derivatives.
Finally, we expect that the formulas and identities obtained in this paper
will find applications in mathematical physics, engineering and related
areas of applied analysis, where hypergeometric–type solutions with tunable
parameters occur naturally.

\section*{Affiliations}
\address{{\bf Ayman Shehata}: Department of Mathematics, Faculty of Science, Assiut University, Assiut 71516, Egypt.}\\
\email{aymanshehata@science.aun.edu.eg, drshehata2009@gmail.com, drshehata2006@yahoo.com}\\
{\bf ORCID ID:  0000-0001-9041-6752}

\address{{\bf Recep \c{S}ahin}: Department of Mathematics,  Faculty of Arts and Sciences, K\i r\i kkale University, 71450, K\i r\i kkale, Turkey.}\\
\email{recepsahin@kku.edu.tr}\\
{\bf ORCID ID: 0000-0001-5713-3830 }

\address{{\bf O\u{g}uz Ya\u{g}c\i}: Department of Mathematics,  Faculty of Arts and Sciences, K\i r\i kkale University, 71450, K\i r\i kkale, Turkey.}\\
\email{oguzyagci26@gmail.com, 1588151031@kku.edu.tr}\\
{\bf ORCID ID: 0000-0001-9902-8094}

\address{{\bf Shimaa I. Moustafa}: Department of Mathematics, Faculty of Science, Assiut University, Assiut 71516, Egypt.}\\
\email{shimaa1362011@yahoo.com, shimaa$_{-}$m.@science.aun.edu.eg}\\
{\bf ORCID ID:  0000-0001-8589-9948}


\begin{thebibliography}{99}
	
	\bibitem{AncaraniGasaneo2008}
	L.~U.~Ancarani and G.~Gasaneo,
	Derivatives of any order of the confluent hypergeometric function ${}_1F_{1}(a;b;z)$ with respect to the parameter $a$ or $b$,
	\emph{J. Math. Phys.} \textbf{49} (2008), 063508.
	
	\bibitem{AncaraniGasaneo2009}
	L.~U.~Ancarani and G.~Gasaneo,
	Derivatives of any order of the Gaussian hypergeometric function ${}_2F_{1}(a;b;c;z)$ with respect to the parameters $a$, $b$ and $c$,
	\emph{J. Phys. A: Math. Theor.} \textbf{42} (2009), 395208.
	
	\bibitem{AncaraniGasaneo2010}
	L.~U.~Ancarani and G.~Gasaneo,
	Derivatives of any order of the hypergeometric function ${}_pF_{q}(a_1,\dots,a_p;b_1,\dots,b_q;z)$ with respect to the parameters $a_i$ and $b_i$,
	\emph{J. Phys. A: Math. Theor.} \textbf{43} (2010), 085210.
	
	\bibitem{AncaraniDelPuntaGasaneo2017}
	L.~U.~Ancarani, J.~A.~Del~Punta and G.~Gasaneo,
	Derivatives of Horn hypergeometric functions with respect to their parameters,
	\emph{J. Math. Phys.} \textbf{58} (2017), 073504.
	
	\bibitem{AppellKampedeFeriet}
	P. Appell, and J. Kamp{\'e} de F{\'e}riet,
	\emph{Fonctions hyperg{\'e}om{\'e}triques et hypersph{\'e}riques: polyn{\^o}mes d'Hermite}, Gauthier-Villars, Paris, 1926.
	
	\bibitem{Fejzullahu2017}
	B.~Xh.~Fejzullahu,
	Parameter derivatives of the generalized hypergeometric function,
	\emph{Integral Transforms Spec. Funct.} \textbf{28}(11) (2017), 781--788.
	
	\bibitem{BytevKniehlMoch2017}
	V.~Bytev, B.~Kniehl and S.~Moch,
	Derivatives of Horn-type hypergeometric functions with respect to their parameters,
	preprint, submitted 18 December 2017.
	
	\bibitem{ExtHyperRefs}
	H.~Exton, \emph{Multiple Hypergeometric Functions and Applications}, Ellis Horwood Limited, Chichester, 1976.
	
	\bibitem{Froehlich1994}
	J.~Froehlich,
	Parameter derivatives of the Jacobi polynomials and Gaussian hypergeometric function,
	\emph{Integral Transforms Spec. Funct.} \textbf{2} (1994), 252--266.
	
	\bibitem{KangAn2015}
	H.~Kang and C.~An,
	Differentiation formulas of some hypergeometric functions with respect to all parameters,
	\emph{Appl. Math. Comput.} \textbf{258} (2015), 454--464.
	
	\bibitem{Sahin2015}
	R.~\c{S}ahin,
	Recursion formulas for Srivastava's hypergeometric functions,
	\emph{Math. Slovaca} \textbf{65}(6) (2015), 1345--1360.
	
	\bibitem{SahinYagci2018}
	R.~\c{S}ahin and O.~Ya\u{g}c{\i},
	$H^{\tau_1,\tau_2,\tau_3}_A$ Srivastava hypergeometric function,
	\emph{Math. Sci. Appl. E-Notes} \textbf{6}(2) (2018), 1--9.
	
	\bibitem{SahatVerma2015}
	V.~Sahat and A.~Verma,
	Derivatives of Appell functions with respect to parameters,
	\emph{J. Inequal. Spec. Funct.} \textbf{6}(3) (2015), 1--16.
	
	\bibitem{SofotasiosBrychkov2018}
	P.~C.~Sofotasios and Yu.~A.~Brychkov,
	On derivatives of hypergeometric functions and classical polynomials with respect to parameters,
	\emph{Integral Transforms Spec. Funct.} \textbf{29}(11) (2018), 852--865.
	
	\bibitem{Srivastava1964}
	H.~M.~Srivastava,
	Hypergeometric functions of three variables,
	\emph{Ganita} \textbf{15}(2) (1964), 97--108.
	
	\bibitem{SrivastavaManocha1984}
	H.~M.~Srivastava and H.~L.~Manocha,
	\emph{A Treatise on Generating Functions},
	Halsted Press (Ellis Horwood Ltd., Chichester),
	John Wiley \& Sons, New York, Chichester, Brisbane and Toronto, 1984.
	
	\bibitem{SrivastavaKarlsson1985}
	H.~M.~Srivastava and P.~W.~Karlsson,
	\emph{Multiple Gaussian Hypergeometric Series},
	Halsted Press (Ellis Horwood Ltd., Chichester),
	John Wiley \& Sons, New York, Chichester, Brisbane and Toronto, 1985.
	
	\bibitem{PathanSrivastava}
	H.~M.~Srivastava, M.~A.~Pathan-Yasmeen,
	Some reduction formulas for multiple hypergeometric series,
	\emph{Rendiconti del Seminario Matematico della Universit{\`a} di Padova} \textbf{82} (1989), 1--7.
	
	\bibitem{Pathan1979}
	M.~A.~Pathan,
	On a transformation of a general hypergeometric series of four variables,
	\emph{Nederl. Akad. Wetensch. Proc. Ser. A} \textbf{82},
	\emph{Indag. Math.} \textbf{41} (1979), 171--175.
	
	
	\bibitem{Yagci2019}
	O.~Ya\u{g}c{\i},
	$H^{\tau_1,\tau_2,\tau_3}_B$ Srivastava hypergeometric function,
	\emph{Math. Sci. Appl. E-Notes} \textbf{7}(2) (2019), 195--204.
	
\end{thebibliography}
\end{document}